\documentclass[3p,preprint]{elsarticle}

\usepackage{amsmath}
\usepackage{amsfonts}
\usepackage{amssymb,latexsym}
\usepackage{amsthm}
\usepackage{graphicx}

\newtheorem{theorem}{Theorem}[section]

\newtheorem{Remark}{Remark}[section]

\usepackage{color}
\usepackage{newlfont}
\usepackage{subfigure}
\usepackage{verbatim}
\usepackage{rotating}
\usepackage{multirow}
\usepackage{units}
\usepackage{bm}
\usepackage[english]{babel}
\usepackage[utf8]{inputenc}
\usepackage{algorithm}
\usepackage[noend]{algpseudocode}
\usepackage{booktabs}
\usepackage{caption}
\usepackage{hyperref}

\journal{}

\begin{document}

\begin{frontmatter}

\title{Cheap and stable quadrature on polyhedral elements}

\author[address-PD]{Alvise Sommariva}

\author[address-PD]{Marco Vianello\corref{corrauthor}}
\ead{marcov@math.unipd.it}

\cortext[corrauthor]{Corresponding author}

\address[address-PD]{Department of Mathematics, University of Padova, Italy}

\begin{abstract}
We discuss a cheap tetrahedra-free approach to the numerical integration of polynomials on polyhedral elements, 
based on hyperinterpolation in a bounding box and Chebyshev moment computation via the divergence theorem. No conditioning issues arise, since no matrix factorization or inversion is needed. The resulting quadrature formula is theoretically stable even in the presence of some negative weights. 
\end{abstract}


\end{frontmatter}

\section{Introduction} 
Within the field of polytopal FEM, which has been costantly growing during the last two decades, one of the main computational bottlenecks, especially for high-order methods in 3D, is the necessity of computing in a fast and stable way the integrals of products of polynomials naturally arising on arbitrary polyhedral elements, avoiding sub-tessellation into tetrahedra.  

Indeed, a specific literature on this topic has been 
emerging, 
where we may quote for example, with no pretence of exahustivity,  \cite{AHP18,BN23,CLS15,CS20,HD23,LvPBD24,MS11,SV24,SDAW14,SSVW17} and the references therein. Among them we focus on  tetrahedra-free moment-based quadratures, which rely on the divergence theorem, together with the availability of efficient quadrature formulas for the polygonal faces and a possible compression in case of high-cardinality; cf., e.g., \cite{SV24,SDAW14,SSVW17}.

The recent paper \cite{LvPBD24} proposed an appealing way to reduce the computational cost, that ultimately corresponds to moment-matching with 
the monomial basis by supporting the quadrature formula on approximate Fekete points of an enclosing box for the polyhedral element, computed by QR factorization with column pivoting of a large Vandermonde matrix (cf. \cite{SV09}). In such a way the relevant 
computations can be made once and for all in a reference box, and the formula weights can be computed by solving relatively small linear systems (whose size is the dimension of the exactness polynomial space). The resulting weights are not all positive, but the sums of their absolute values remains experimentally bounded, ensuring empirical stability. 

In the present paper we adopt a similar idea, namely that of supporting the quadrature formula on a low cardinality set, but this time such a set is the support itself of a suitable hyperinterpolation formula in the enclosing box, 
by a sort of generalized  Clenshaw-Curtis approach proposed in \cite{SVZ08}. We recall that hyperinterpolation, introduced by Sloan in the seminal paper \cite{S95}, is a Fourier-like orthogonal projection on a total-degree polynomial space with respect to an absolutely continuous measure, 
discretized by an algebraic quadrature formula with positive weights. 

This approach gives not only a simple explicit and cheap formula for the computation of the quadrature weights, consisting essentially of a single matrix-by-vector product, but also theoretically ensures stability of the quadrature formula. In practice, the computational cost substantially reduces to the mere cost of computing the moments of the orthogonal basis used for hyperinterpolation. Moreover, differently from \cite{LvPBD24}, no conditioning issue can arise by increasing the exactness degree, not only because we can use in a natural way for example the Chebyshev basis instead of the monomial basis, but mainly because no QR factorization or linear system is involved. The theoretical foundation upon hyperinterpolation 
following \cite{SVZ08}, with the consequent stability result, make the present approach also quite different from the adaptive scheme recently proposed in \cite{BN23}.

The paper is organized as follows. In Section 2 we describe the theoretical base of the method, while in Section 3 we discuss its implementation and provide some numerical examples.

\section{Polyhedral quadrature by hyperinterpolation}
The proposed quadrature formula relies on two main results, that were proved in \cite{S95} and \cite {SVZ08}. We summarize for convenience these results as a single theorem, which concerns product-like formulas obtained via hyperinterpolation, and their stability. 

Below, we shall denote by $\mathbb{P}_n^d$ the space of $d$-variate polynomials with total degree not exceeding $n$, with dimension $N={n+d\choose n}$.

\begin{theorem}\label{2.1}
Let $K\subset \mathbb{R}^d$ be a compact subset, $\mu$ an absolutely continuous measure on $K$ with respect to the Lebesgue measure. Denote by $\{\phi_j\}_{1\leq j\leq N}$ an orthonormal polynomial basis of $\mathbb{P}_n^d$ for $\mu$. Moreover, let $(X,\mathbf{u})=(\{P_i\},\{u_i)\})$, $1\leq i\leq \nu=\nu(n)$, be the nodes and positive weights of a quadrature formula for integration in $d\mu$, exact on $\mathbb{P}_{2n}^d$ (the polynomials with total degree not exceeding $2n$), and $h\in L^2_{\mu}(K)$.   

Then, the following algebraic product-like formula holds
\begin{equation}
\int_K{h(P)f(P)\,d\mu} =\sum_{i=1}^\nu{w_i\,f(P_i)}\;,\;\forall f \in \mathbb{P}_n^d\;,
\end{equation}
where the quadrature weights $\{w_i\}$ are defined by the product-like moments
\begin{equation}
w_i=u_i\,\sum_{j=1}^N{\phi_j(P_i)\,m_j}\;,1\leq i\leq \nu\;,\;\;\;m_j=\int_K{\phi_j(P)\,h(P)\,d\mu}\;.
\end{equation}
Moreover, the formula is stable, since
\begin{equation} \label{stable}
\lim_{n\to \infty}\sum_{i=1}^\nu{|w_i|}=\int_K{|h(P)|\,d\mu}\,.
\end{equation}
\end{theorem}

For the proof of this theorem, we refer the reader to \cite{S95} and \cite{SVZ08}. We recall that $\nu \geq N$, and that when $\nu=N$ the quadrature formula for $\mu$ is called minimal; cf. \cite{M75}. Concerning stability, though there could be some negative weights, (\ref{stable}) clearly implies that $\sum_{i=1}^\nu{|w_i|}$ is bounded. Moreover, notice that if $h$ is almost everywhere nonnegative, then the stability parameter of the quadrature formula, namely $\sum_i{|w_i|}/|\sum_i{w_i}|=\sum_i{|w_i|}/\|h\|_{L^1(K)}$, tends to 1 as $n\to \infty$.

We turn now to the main goal of the present paper, that is constructing a cheap and stable quadrature formula for 
\begin{equation} \label{polyhedron}
\int_\Omega{f(P)\,dP}\;,\;\;\Omega\subset \mathbb{R}^3\;\;\mbox{polyhedron}\;,\;\;f\;\mbox{polynomial}\;.
\end{equation}
Here and below, $P=(x,y,z)$ and $dP=dx\,dy\,dz$. 

Now, take $K=B\supseteq \Omega$, where $B$ is a Cartesian bounding box 
for a polyhedron $\Omega$, that up to an affine change of variables can be taken as the cube $B=[-1,1]^3$. Given any absolutely continuous measure $d\mu=\sigma(P)\,dP$, $\sigma\in L^1_+(B)$, for which we know an algebraic quadrature formula with positive weights for total degree $2n$, we can apply Theorem 1 by writing 
\begin{equation} \label{h}
\int_\Omega{f(P)\,dP}=\int_B{h(P)f(P)\sigma(P)\,dP}\;,\;h(P)=I_\Omega(P)/\sigma(P)\;,
\end{equation}
provided that $h\in L^2_{\mu}(B)$, that is $1/\sigma\in L^1_+(\Omega)$. We then obtain
\begin{equation} \label{qformula}
\int_\Omega{f(P)\,dP}=\sum_{i=1}^\nu{w_i\,f(P_i)}\;,\forall f\in \mathbb{P}_n^3\;,\;\;
w_i=u_i\,\sum_{j=1}^N{\phi_j(P_i)\,m_j}\;,1\leq i\leq \nu\;,\;\;\;m_j=\int_\Omega{\phi_j(P)\,dP}\;,
\end{equation}
where $N=dim(\mathbb{P}_n^3)=(n+1)(n+2)(n+3)/6$.
Moreover, in this case we have by (\ref{stable})
\begin{equation} \label{stable2}
 \lim_{n\to \infty}\sum_{i=1}^\nu{|w_i|}=   
 \int_B{(I_\Omega(P)/\sigma(P))\,\sigma(P)\,dP}=vol(\Omega)\;.
\end{equation}

Observe that this approach to quadrature on polyhedra is quite general, and clearly extendable to any compact set where one is able to compute the Lebesgue measure moments $\{m_j\}$ in (\ref{qformula}), for a polynomial basis orthogonal with respect to an absolutely continuous measure in a bounding box. In the case of polyhedra, as we shall see below, 
this can be effectively done by the divergence theorem. 

We have several natural choices for the measure $\mu$. The first that comes to mind is simply $d\mu=dP$, the Lebesgue measure itself. In such a way, the orthogonal basis is the total-degree product Legendre basis (cf. e.g. \cite{DX01}), and as quadrature formula 
of exactness degree $2n$ we can choose the tensorial Gauss-Legendre rule  
with $\nu=(n+1)^3$ nodes, or even a minimal ($\nu=N$) or near-minimal formula for the degrees where it is available in the quadrature literature (cf. e.g. \cite{C03}). 

It is worth stressing however that the choice of the underlying quadrature formula in the cube is not decisive, because the cardinalities are in any case small enough to make largely predominant the cost of moment computation, via the divergence theorem and algebraic quadrature on the polygonal faces of the polyhedral surface (as we shall see in the numerical section). On the other hand, the computation of the necessary primitives is facilitated by known analytical formulas for the classical univariate orthogonal polynomials.

But we can also choose the product Chebyshev measure. In this case the orthogonal basis is the total-degree product Chebyshev basis, and the tensorial Gauss-Chebyshev 
or Gauss-Chebyshev-Lobatto rules with $\nu=(n+1)^3$ nodes become a natural choice. A lower cardinality formula could be adopted, for example the formula proposed in \cite{DMVX09} with $\nu\approx (n+1)^3/4$, or even a minimal or near-minimal formula at the available degrees in the quadrature literature (cf. e.g. \cite{PV15}), but again this is not really relevant for a significative reduction of the overall cost.    

\begin{Remark}
We stress that the present approach, like that developed in \cite{LvPBD24}, gives not only a numerical integration method for polynomials, but also a quadrature formula on a polyhedron $\Omega$, applicable 
to any function that is defined and well-approximated by polynomials in the whole bounding box $B\supseteq \Omega$. Standard estimates in quadrature theory together with (\ref{stable2}) allow indeed to write the following error bound for every $f\in C(B)$
\begin{equation} \label{errest}
\left|\int_\Omega{f(P)dP}-\sum_{i=1}^\nu{w_i\,f(P_i)}\right|\leq \left(vol(\Omega)+\sum_{i=1}^\nu{|w_i|}\right)\,E_n(f;B)\sim 2\,vol(\Omega)\,E_n(f;B)\;,
\end{equation}
where $E_n(f;B)=\inf_{\phi\in \mathbb{P}_n}\max_{P\in B}{|f(P)-\phi(P)|}$. The decay rate of 
$E_n(f;B)$ as $n\to \infty$ depends on the 
regularity of $f$, by a multivariate version of Jackson theorem; in particular, if $f\in C^{k+1}(B)$ then $E_n(f;B)=\mathcal{O}(n^{-k})$, cf. e.g. \cite{P09}. From this point of view, we mark a difference with respect to other methods like \cite{AHP18},  
that are by construction restricted to polynomials.
\end{Remark}

\section{Implementation and numerical tests} 
The computational steps necessary for the construction of the quadrature formula (\ref{qformula}), exact for $\mathbb{P}_n^3$ on an 
arbitrary polyhedron $\Omega$, can be briefly summarized as follows: 

\begin{itemize}
\item[$(i)$] determine a Cartesian bounding box for the polyhedron and 
compute the nodes $\{P_i\}$ and weights $\mathbf{u}=\{u_i\}$, $1\leq i\leq \nu$, of a quadrature formula exact for 
$\mathbb{P}_{2n}^3$ for a given absolutely continuous measure 
$d\mu=\sigma(P)dP$ on the bounding box; 

\item[$(ii)$] compute an orthonormal basis $\{\phi_1,\dots,\phi_N\}$ of $\mathbb{P}_n^3$ with respect to $d\mu$, and the corresponding 
Lebesgue moments $\mathbf{m}=\{m_1,\dots,m_N\}$, $m_j=\int_\Omega{\phi_j(P)\,dP}$, by the divergence theorem via the (oriented) planar polygonal faces of the polyhedron $$m_j=\int_{\partial \Omega}{\phi_j(P)\,n_1(P)\,dS}=\sum_{faces}\int_{face}{\phi_j(P)\,n_1^{face}\,dS}\;,$$
where $dS$ is the surface measure, $\phi_j(P)=\int{\phi_j(P)\,dx}$ is a primitive with respect for example to $x$, that is $\partial_x \phi_j(P)=\phi_j(P)$, and $n_1$ is the first component of the outward normal vector to the polyhedral surface, which is constant 
on each planar face;

\item[$(iii)$] form the Vandermonde-like matrix $V=V_n(\{P_i\})=[\phi_j(P_i)]\in \mathbb{R}^{\nu\times N}$ and compute the final weights as 
a scaling of a matrix-by-vector product as 
$$
\mathbf{w}=diag(\mathbf{u})\,V\mathbf{m}\;,
$$
or in a Matlab-like notation $\mathbf{w}=\mathbf{u}.*\,V\mathbf{m}\;$.

\end{itemize}

Notice that $(i)$ can be made completely independent of the polyhedron, 
by choosing a reference box such as $[-1,1]^3$, via an affine change of variables (namely, a translation plus scaling of variables) which affects the integrals and the weights only by multiplicative constants. On the other hand, $(iii)$ depends on the polyhedron only via the moment vector $\mathbf{m}$ computed in $(ii)$, because 
\begin{itemize}
\item {\em the Vandermonde-like matrix can be computed once and for all in the reference cube}.  
\end{itemize}
This aspect is particularly relevant in the application to polyhedral FEM, where quadrature has to be applied to a potentially very large number of different polyhedral elements.

Concerning $(ii)$, assuming again with no loss of generality that $B=[-1,1]^3$, the computation is substantially simplified by choosing 
a product type orthogonal basis, like e.g. the Chebyshev basis $$\phi_j(P)=c_j\,T_{\alpha}(x)T_{\beta}(y)T_{\gamma}(z)$$ where $(\alpha,\beta,\gamma)\in \{0,1,\dots,n\}^3$, $0\leq \alpha+\beta+\gamma\leq n$, the $c_j$ are normalization constants, and the index $j$ corresponds to a suitable ordering of the triples, for example a lexicographical ordering. 
In particular, if $j$ corresponds to a certain triple $(\alpha,\beta,\gamma)$, then $c_j=a_\alpha a_\beta a_\gamma$, where $a_0=1/\sqrt{\pi}$ and $a_k=\sqrt{2/\pi}$ for $k>0$. Indeed, a primitive in this case is analytically known, 
$$\phi_j(P)=c_j\,\int_0^x{T_{\alpha}(x)dx}\;T_{\beta}(y)T_{\gamma}(z)$$
where $\int_0^x{T_{\alpha}(x)dx}=\frac{T_{\alpha+1}(x)}{2(\alpha+1)}
-\frac{T_{\alpha-1}(x)}{2(\alpha-1)}$ for $\alpha\geq 2$, whereas 
trivially $\int_0^x{T_{0}(x)dx}=x$, $\int_0^x{T_{1}(x)dx}=x^2/2$; cf. e.g. \cite{MH02}. 
The integrals on the planar faces in $(ii)$ can be computed in several ways, for example by the formulas based on piecewise product Gauss-Legendre quadrature, developed in \cite{SV07}. For a complete discussion on the computation of the polyhedral moments for the product Chebyshev basis, we refer the reader to \cite[\S 2.1]{SV24}.

Finally, the cost of $(iii)$ is substantially that of a matrix-by-vector product. 
Even though it is not the dominant one, with the choice of the total-degree product Chebyshev basis and of the tensorial Gauss-Chebyshev rule, this cost could be even lowered since it essentially 
corresponds to a discrete cosine transform of the moment vector, that can be accelerated by the FFT. 

The numerical tests below have been performed using Matlab R2024A on an Intel Core Ultra 5 125H processor, with frequency 3.60 GHz and 16 GB of RAM. A preliminary non optimized version of the Matlab code, named {\em cheapQ}, is available at \cite{software}. 
In our implementation, to compute a {\it{cheap}} formula with degree of exactness $n$ on $\Omega$, in accordance to Theorem 2.1 we used as set of nodes those of a tensorial Gauss-Chebyshev rule on $[-1,1]^3$ with degree of exactness $2n$, namely $(n+1)^3$ nodes, scaled to the bounding box.

\begin{figure}[!htbp]
  \centering
   {\includegraphics[scale=0.40,clip]{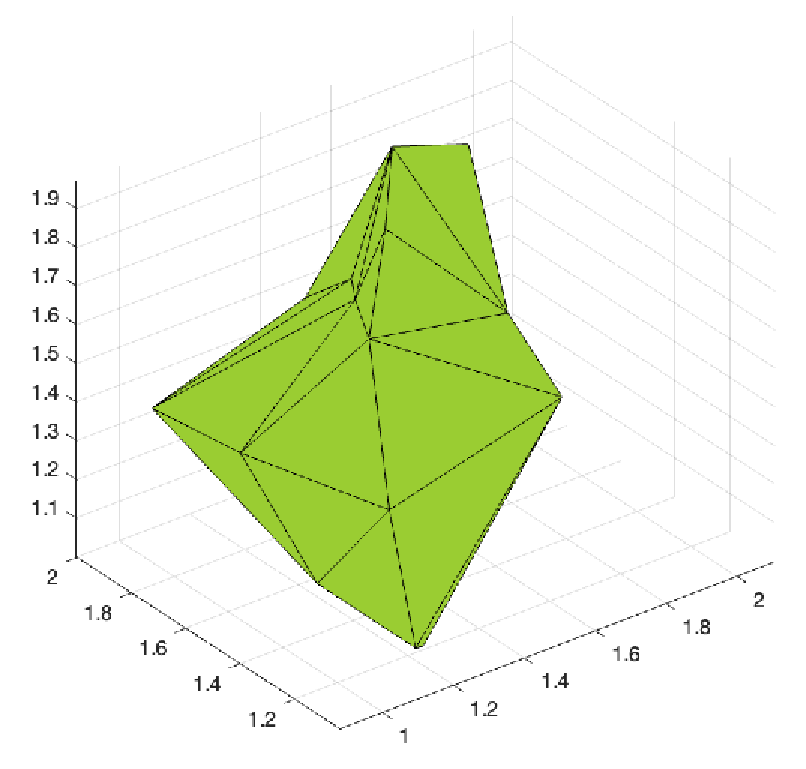}}
   {\includegraphics[scale=0.40,clip]{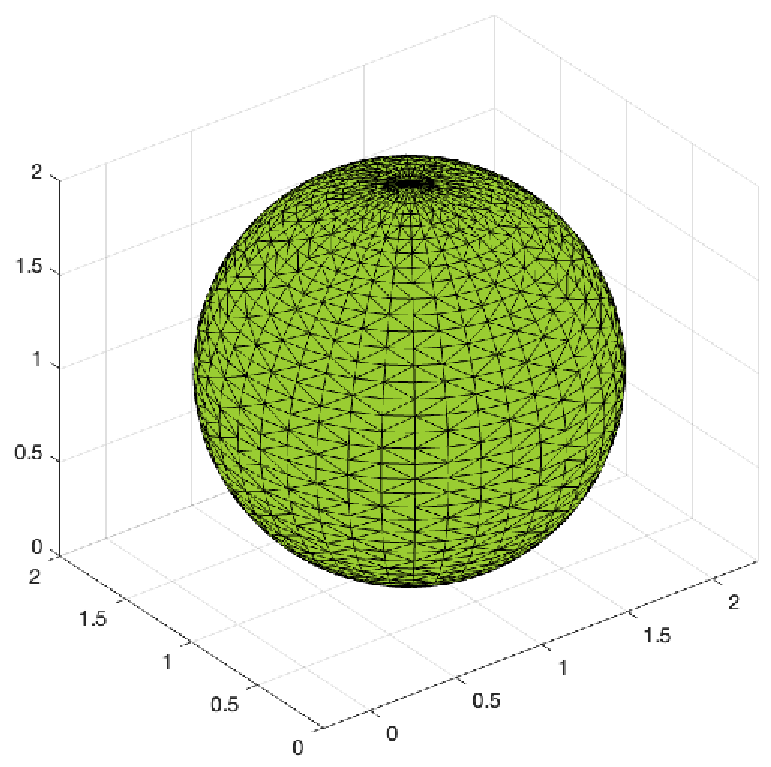}}
      {\includegraphics[scale=0.40 ,clip]{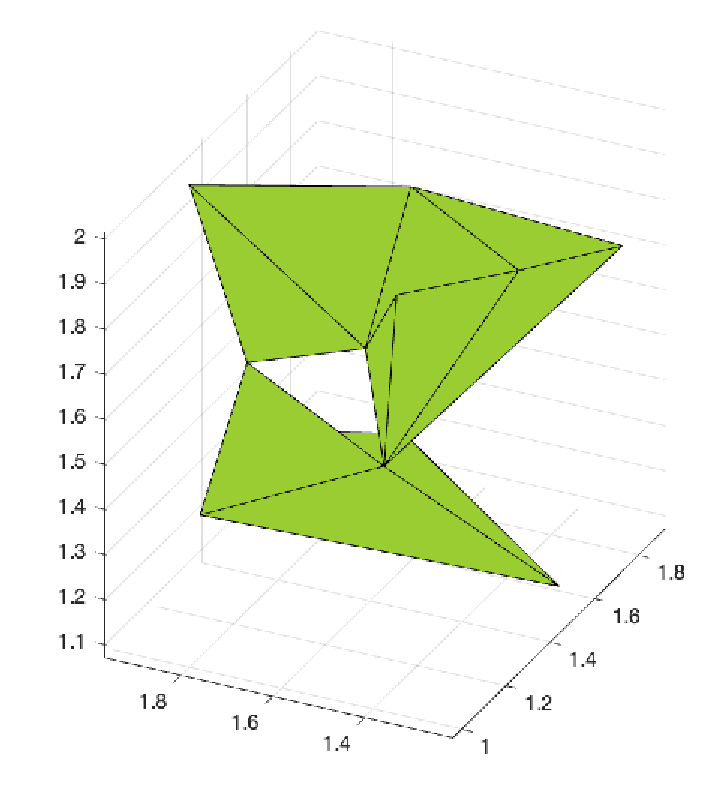}}
 \caption{Examples of polyhedral domains. {{Left:}} $\Omega_1$ (nonconvex, 20 facets); {{Center}}: $\Omega_2$ (convex, 760 facets); {{Right}}: $\Omega_3$ (multiply connected, 20 facets).}
 \label{fig_2CRI}
 \end{figure}

\begin{figure}[!htbp]
  \centering
   {\includegraphics[scale=0.40,clip]{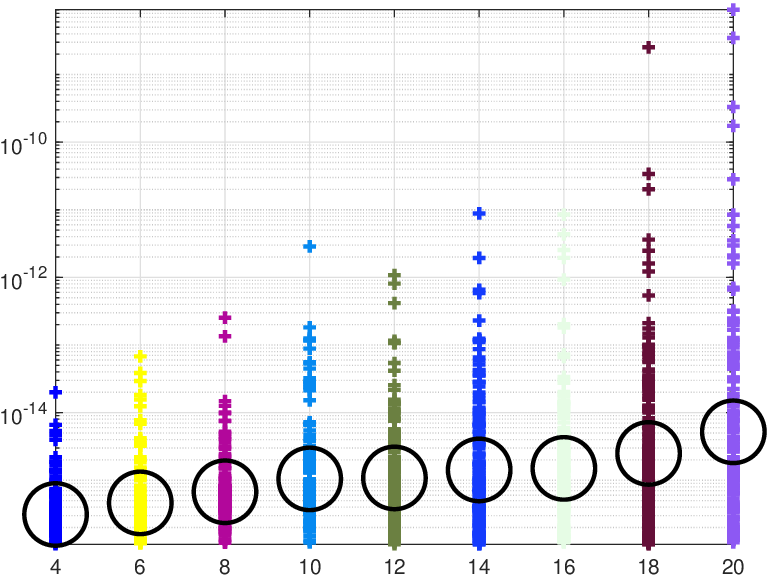}}
   {\includegraphics[scale=0.40,clip]{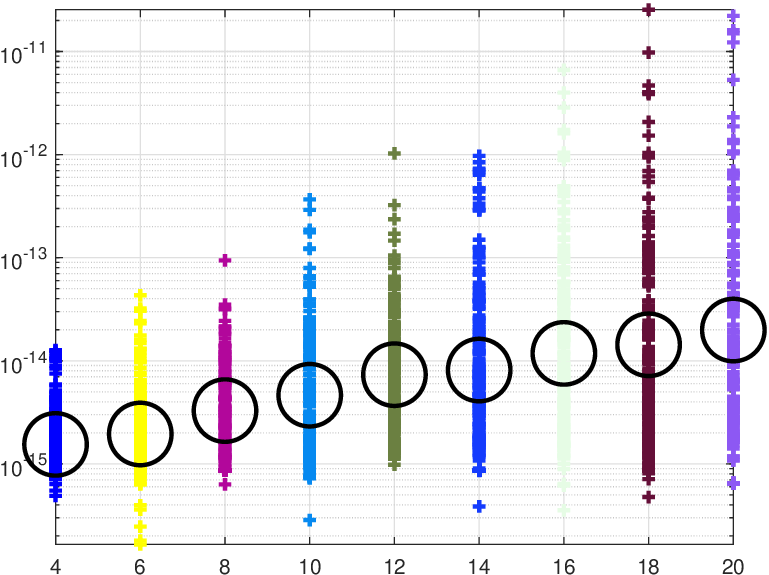}}
      {\includegraphics[scale=0.40,clip]{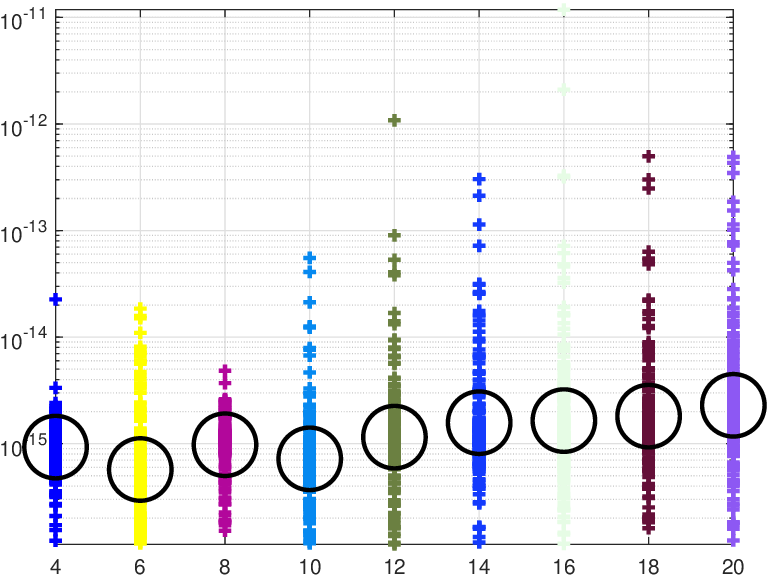}}

\caption{Relative errors $E(g_k)$ of the tetrahedra-free rule over 200 polynomial integrands of the form $g_k=(a_kx+b_ky+c_kz+d_k)^n$ on the three polyhedra of Figure 1, where $a_k,b_k,c_k,d_k$ are uniform random coefficients in $[-1,1]$ and $n=4,6,8,\dots,20$; the circles correspond to the average logarithmic error $\sum_{k=1}^{200}{\log(E(g_k))/200}$.}
 \label{fig_3CRI}
 \end{figure}

In Figure {\ref{fig_3CRI}} we show the relative integration errors 
for the random polynomials $g_k(x,y,z)=(a_kx+b_ky+c_kz+d_k)^n$ on the three polyhedral domains of Figure \ref{fig_2CRI}, where $a_k,b_k,c_k,d_k$ are uniform random coefficients in $[-1,1]$. The domains have been obtained starting from suitable point clouds by the Matlab built-in command {\tt{alphashape}}, which also provides the boundary facets. 
 For each even $n$ in the range between $4$ and $20$, we have made 200 tests and computed the average logarithmic relative error $\sum_{k=1}^{200}{\log(E(g_k))/200}$, displayed by a black circle. The relative errors $E(g_k)$ have been computed using as reference integral the value produced by a tessellation-based formula with exactness degree $n$, taken from \cite{SV24}. We have chosen even degree of exactness $n$, since in a possible application to FEM these correspond to polynomial elements of degree $n/2$. The numerical experiments show that as expected, such average errors are not far from machine precision.

Next in Table 1 we list the average cputimes in seconds of the {\em CheapQ} algorithm. We notice that, as expected, the results for $\Omega_2$ are higher in view of the larger number of facets of the domain. On the other hand, with the two polyhedral elements with 20 facets, we see a good performance, with times ranging from the order of $10^{-3}$ seconds per element for the lowest degrees to $10^{-2}$ seconds for degree 20. 
These times are up to one order of magnitude lower than those of the tetrahedra-free Tchakaloff-like formulas with positive weights and interior nodes, that we have recently constructed in \cite{SV24} by moment-matching with NonNegative Least Squares. 
Finally in Table 2 we report the stability ratios  $\sum_{j=1}^{\nu} |w_j|/vol(\Omega_i)$. We can see that these quantities do not exceed 2 and tend to decrease towards 1 as $n$ increases, as suggested by Theorem {\ref{2.1}}.

\begin{table}[ht]
\begin{center}\footnotesize
{\renewcommand{\arraystretch}{1.1}
\begin{tabular}{| c | c c c c c c c c c | }
\hline
{\mbox{deg} $n$} & $4$ & $6$ & $8$ & $10$ & $12$ & $14$ & $16$ & $18$ & $20$ \\
\hline
$ \Omega_1 $ & 1.2-e03 & 1.4-e03 & 1.7e-03 & 2.3-e03 & 3.4e-03 & 5.1e-03 & 7.7e-03 & 1.9e-02 & 3.4e-02 \\
$ \Omega_2 $ & 3.0e-02 & 3.4e-02 & 4.3e-02 & 5.9e-02 & 8.2e-02 & 1.2e-01 & 1.8e-01 & 4.4e-01 & 9.7e-01 \\
$ \Omega_3 $ & 8.1e-04 & 9.0e-04 & 1.1e-03 & 1.7e-03 & 2.3e-03 & 3.5e-03 & 5.4e-03 & 1.3e-02 & 2.6e-02 \\
 \hline
\end{tabular}
        }
\caption{Average cputimes of {\em CheapQ} on the domains of Fig. 1, varying the algebraic degree of exactness.}
\label{tab1}
\end{center}
\end{table}

\begin{table}[ht]
\begin{center}\footnotesize
{\renewcommand{\arraystretch}{1.1}
\begin{tabular}{| c | c c c c c c c c c | }
\hline
{\mbox{deg} $n$} & $4$ & $6$ & $8$ & $10$ & $12$ & $14$ & $16$ & $18$ & $20$ \\
\hline
$ \Omega_1 $ & 1.55 & 1.40 & 1.30 & 1.25 & 1.23 & 1.21 & 1.19 & 1.17 & 1.17 \\
$ \Omega_2 $ & 1.30 & 1.14 & 1.21 & 1.12 & 1.13 & 1.12 & 1.10 & 1.10 & 1.09 \\
$ \Omega_3 $ & 1.63 & 1.81 & 1.89 & 1.86 & 1.82 & 1.79 & 1.74 & 1.67 & 1.63 \\
 \hline
\end{tabular}
        }
\caption{The ratios $\sum_{j=1}^{\nu} |w_j|/vol(\Omega_i)$ for {\em CheapQ} on the domains of Fig. 1, varying the algebraic degree of exactness.}
\label{tab1}
\end{center}
\end{table}

\section{Conclusions}
We have implemented a quadrature formula without sub-tessellation, which is exact for polynomials up to a given degree on polyhedral elements. The formula is 
based on hyperinterpolation in a bounding box and Chebyshev moment computation via the divergence theorem. The computational bulk is given by computation of the Chebyshev moments, since the final weights are obtained via a matrix-by-vector product where the matrix is element-independent, and can be computed once and for all. No conditioning issues arise, since no matrix factorization or inversion is needed. Moreover, the resulting quadrature formula is theoretically stable even in the presence of some negative weights. 
We are confident that the present method  could become an additional tool, with some improved features with respect to \cite{LvPBD24} and other techniques adopted in the literature, for the efficient and stable computation of stiffness and mass matrices within polyhedral Finite Elements. 
Cheap but still accurate assembly of such matrices, due to the increasing  adoption of polytopal FEM simulations, could have a non negligible fall-out on large scale numerical modelling.

\vskip0.5cm 
\noindent
{\bf Acknowledgements.} 

Work partially supported by the DOR funds of the University of Padova and by the INdAM-GNCS 2024 Project “Kernel and polynomial methods for approximation and integration: theory and application software''.
This research has been accomplished within the Community of Practice 
``Green Computing" of the Arqus European University Alliance, the RITA ``Research ITalian network on Approximation", the SIMAI Activity Group ANA\&A, and the UMI Group TAA ``Approximation Theory and Applications".


\begin{thebibliography}{99}


\bibitem{AHP18} 
P.F. Antonietti, P. Houston, G. Pennesi, 
Fast numerical integration on polytopic meshes with applications to discontinuous Galerkin finite element methods, 
J. Sci. Comput. 77 (2018) 339--1370.

\bibitem{BN23} B. Boroomand, N. Niknejadi, Adaptive quadrature/cubature rule: Application to polytopes, Comp. Methods Appl. Mech. Engrg. 403 (2023), 115726.

\bibitem{CGH14} A. Cangiani, E.H. Georgoulis, P. Houston, 
hp-version discontinuous Galerkin methods on polygonal and polyhedral meshes, Math Models Methods Appl. Sci. 24 (2014), 2009--2041.

\bibitem{CLS15} E.B. Chin, J.B. Lasserre, N.Sukumar, Numerical integration of homogeneous functions on convex and nonconvex polygons and polyhedra, Comput. Mech. 56 (2015), 967--981.

\bibitem{CS20} E.B. Chin, N. Sukumar, An efficient method to integrate polynomials over polytopes and curved solids, Comput. Aided Geom. Des. 82 (2020), 101914. 

\bibitem{C03} R. Cools, An Encyclopaedia of Cubature Formulas, J. Complexity  19 (2003), 445--453 (\url{https://nines.cs.kuleuven.be/research/ecf/}).

\bibitem{DMVX09} S. De Marchi, M. Vianello, Y. Xu, New cubature formulae and hyperinterpolation in three variables, BIT Numer. Math. 49 (2009), 55--73. 

\bibitem{DX01} C.F. Dunkl and Y. Xu, Orthogonal Polynomials of Several Variables,
Encyclopedia of Mathematics and its Applications, vol. 81, Cambridge
University Press, Cambridge, 2001.

\bibitem{HD23} S. Hubrich, A. D\"{u}ster, Numerical integration for nonlinear problems of the finite cell method using an adaptive scheme based on
moment fitting, Comput. Math. Appl. 77 (2019), 1983--1997.

\bibitem{LvPBD24} C. Langlois, T. van Putten, H. B\'{e}riot, E. Deckers, Frugal numerical integration scheme for polytopal domains, Eng. Comput. (2024). 

\bibitem{MH02} J.C. Mason, D.C. Handscombe, Chebyshev Polynomials, Chapman\&Hall/CRC, 2002. 

\bibitem{MW} 
The Mathworks, {\tt{alphaShape}}, 
Polygons and polyhedra from points in 2-D and 3-D, 
{\url{https://www.mathworks.com/help/matlab/ref/alphashape.html}}.

\bibitem{M75} H.M. M\"{o}ller, Kubaturformeln mit minimaler Knotenzahl, Numer. Math.
25 (1975/76), no. 2, 185–200.

\bibitem{MS11} 
S.E. Mousavi, N. Sukumar, Numerical integration of polynomials and discontinuous functions on irregular convex polygons and polyhedrons, 
Comput. Mech. 47 (2011) 535--554.

\bibitem{P09} W. Ple\'{s}niak, Multivariate Jackson Inequality, J. Comput. Appl. Math. 233 (2009), 815--820.

\bibitem{PV15} D. Potts, T. Volkmer, Fast and exact reconstruction of arbitrary multivariate algebraic polynomials in Chebyshev form, 2015 International Conference on Sampling Theory and Applications (SampTA), IEEE, pp. 392--396.

\bibitem{S95} I.H. Sloan, Interpolation and Hyperinterpolation over General regions, J. Approx. Theory 83 (1995) 238--254.

\bibitem{software} A. Sommariva,  {\em CheapQ}, Matlab codes for cheap and stable tetrahedra-free quadrature on polyhedral elements, {\url{https://github.com/alvisesommariva/CheapQ}}.

\bibitem{SV07} A. Sommariva, M. Vianello, Product Gauss cubature over polygons based on Green's integration formula, 
BIT Numerical Mathematics 47 (2007), 441--453.

\bibitem{SV09} A. Sommariva, M. Vianello, Computing approximate Fekete points by QR factorizations of Vandermonde matrices, 
Comput. Math. Appl. 57 (2009), 1324--1336.

\bibitem{SV24} A. Sommariva, M. Vianello, TetraFreeQ: tetrahedra-free quadrature on polyhedral elements, Appl. Numer. Math. 200 (2024), 389--398.

\bibitem{SVZ08} A. Sommariva, M. Vianello, R. Zanovello, Nontensorial Clenshaw-Curtis cubature, 
Numer. Algorithms 49 (2008), 409--427.

\bibitem{SDAW14} Y.
Sudhakar, J.P. Moitinho De Almeida, W.A. Wall, An accurate, robust, and easy-to-implement method for integration over arbitrary polyhedra: application to embedded interface methods, J. Comput. Phys. 273 (2014), 393--415. 

\bibitem{SSVW17} 
Y. Sudhakar, A. Sommariva, M. Vianello, W.A. Wall, 
On the use of compressed polyhedral quadrature formulas in embedded interface methods, SIAM J. Sci. Comput. 39 (2017) B571--B587.


\end{thebibliography}
\end{document}